\documentclass{article}%
\usepackage{amssymb}
\usepackage{amsfonts}
\usepackage{amsmath}%
\usepackage{graphicx}
\providecommand{\U}[1]{\protect\rule{.1in}{.1in}}

\pdfoutput=1
\begin{document}

\title{Analytical Representations\\of Divisors of Integers}
\author{Krzysztof Ma\'{s}lanka\\e-mail: krzysiek2357@gmail.com\bigskip\\
Institute for the History of Science\\Polish Academy of Sciences}
\maketitle

\begin{abstract}
Certain analytical expressions which "feel" the divisors of natural numbers
are investigated. We show that these expressions encode to some extent the
well-known algorithm of the sieve of Eratosthenes.

Most part of the text is written in pedagogical style, however some formulas
are new.\bigskip

\bigskip\bigskip

\end{abstract}

\section{Notation and Conventions}

Throughout this paper we shall adopt the following notation and conventions:
$n$ is a given natural number and $k$ is a possible divisor of $n$. If $k$
actually divides $n$ then $j=n/k$. Let $f(x)$ denotes any real analytic
function defined in the neighborhood of the origin by a power series%
\begin{equation}
f(x)=\sum_{j=0}^{\infty}c_{j}x^{j}\label{f0}%
\end{equation}
with all $c_{j}\neq0$ ($i=1,2,3...$). It will be shown that $j$ is also the
exponent of $x$ in the expansion (\ref{f0}) around zero and $j$ labels
half-lines or rays of divisors (see below).

\section{Motivation}

The theory of divisors of integers is the cornerstone of elementary number
theory. It is convenient to introduce \textit{the characteristic function for
divisors}:\bigskip

\textbf{Definition.} For any $n,k\in%
\mathbb{N}
$%
\begin{equation}
\overset{\_}{\alpha}_{nk}:=\left\{
\genfrac{}{}{0pt}{}{1\quad\text{if}\quad k\mid n}{0\quad\text{if}\quad k\nmid
n}%
\right. \label{alpha1}%
\end{equation}

Another pretty obvious (and rather useless in numerical calculations)
representation of (\ref{alpha1}) is:%
\begin{equation}
\overset{\_}{\alpha}_{nk}=\frac{1}{\Gamma\left(  1-\operatorname{mod}\left(
n,k\right)  \right)  }\label{alpha2}%
\end{equation}
where $\Gamma(s)$\ denotes the Euler gamma function and $\operatorname{mod}%
\left(  n,k\right)  $ gives the remainder on division of $n$ by $k$. In fact
(\ref{alpha2}) is more general than (\ref{alpha1})\ since it may be calculated
also for non-integer or even complex values of $n$ and $k$ but this leads to
some interpretation difficulties which we shall not discuss here.

Consider the following expression for some natural numbers $n$ and $k$:%
\begin{equation}
\alpha_{nk}=\left.  \frac{\mathrm{d}^{n}}{\mathrm{d}x^{n}}e^{x^{k}}\right\vert
_{x=0}\label{alpha3}%
\end{equation}
We will prove the following\bigskip

\textbf{Theorem.} Apart from a trivial normalization factor,
$\overset{\_}{\alpha}_{nk}$\ defined in formula\ (\ref{alpha1}) is equal to
$\alpha_{nk}$ defined in\ (\ref{alpha3}).\bigskip

\textbf{Proof.} Expanding the exponential function in (\ref{alpha3})\ in power
series and performing term-by-term differentiation we get:%
\begin{equation}
\alpha_{nk}=\left.  \frac{\mathrm{d}^{n}}{\mathrm{d}x^{n}}%
{\displaystyle\sum\limits_{j=0}^{\infty}}
\frac{\left(  x^{k}\right)  ^{j}}{j!}\right\vert _{x=0}=\left.
{\displaystyle\sum\limits_{j=0}^{\infty}}
\frac{1}{j!}\frac{\mathrm{d}^{n}}{\mathrm{d}x^{n}}x^{jk}\right\vert
_{x=0}\label{alpha4}%
\end{equation}
Recall the general formulas for the $n$-th\ derivative of $x^{p}$ with respect
to $x$%
\begin{equation}
\frac{\mathrm{d}^{n}}{\mathrm{d}x^{n}}x^{p}=\frac{\Gamma\left(  p+1\right)
}{\Gamma\left(  p+1-n\right)  }x^{p-n}=n!\binom{p}{n}x^{p-n}%
\label{power derivative 1}%
\end{equation}%
\begin{equation}
\frac{\mathrm{d}^{n}}{\mathrm{d}x^{n}}x^{p}=\left(  -1\right)  ^{n}%
\frac{\Gamma\left(  n-p\right)  }{\Gamma\left(  -p\right)  }x^{p-n}%
\label{power derivative 2}%
\end{equation}
where the second formula stems from properties of the gamma function and is
suitable for integer negative $p$ (see e.g. \cite{MillerRoss}). Note that the
order of derivative $n$ does not have to be integer but for integer $n$\ both
(\ref{power derivative 1}) and (\ref{power derivative 2}) reduce to the
well-known elementary differentiation rule. Using (\ref{power derivative 1})
we get:%
\begin{equation}
\alpha_{nk}=\left.
{\displaystyle\sum\limits_{j=0}^{\infty}}
\frac{1}{j!}\frac{\Gamma\left(  jk+1\right)  }{\Gamma\left(  jk+1-n\right)
}x^{jk-n}\right\vert _{x=0}\left.  =n!%
{\displaystyle\sum\limits_{j=0}^{\infty}}
\frac{1}{j!}\binom{jk}{n}x^{jk-n}\right\vert _{x=0}\label{alpha5}%
\end{equation}
By simple inspection of (\ref{alpha5}) we see why this expression "feels" the
divisors of the integer $n$. Indeed, when taking the limit $x\rightarrow0$ the
only non-zero term in the series appears when $jk=n$ for some integer $j$, and
this occurs if and only if $k$ divides $n$. All terms with $jk>n$ disappear in
the limit $x\rightarrow0$ whereas those with $jk<n$, although singular in
$x=0$, vanish since the binomial coefficient term is zero. Therefore, in the
summation (\ref{alpha5}) at most only one term can survive in the limit
process.$\blacksquare$

The above reasoning might appear far too excessive. However, it guarantees
that among divisors none have been omitted. It should also be stressed that it
may be used as a starting point for various generalizations since $n$ need not
to be integer.\bigskip

It is easy to guess the normalizing factor:%
\begin{equation}
\alpha_{nk}=\frac{1}{n!}\left(  \frac{n}{k}\right)  !\left.  \frac
{\mathrm{d}^{n}}{\mathrm{d}x^{n}}e^{x^{k}}\right\vert _{x=0}\label{alpha6}%
\end{equation}
Using the same reasoning we can derive similar expression for $\alpha_{nk}$:%
\begin{equation}
\alpha_{nk}=\frac{\left(  k!\right)  ^{n/k}}{n!}\left(  \frac{n}{k}\right)
!\left.  \frac{\mathrm{d}^{n}}{\mathrm{d}x^{n}}e^{\frac{x^{k}}{k!}}\right\vert
_{x=0}\label{alpha7}%
\end{equation}

\section{Simple example}

In a natural way coefficients $\alpha_{nk}$\ may be regarded as a square
matrix of arbitrarily large dimension where the running integer $n$ labels
rows and the potential divisor $k$ labels columns. The entries of this matrix
are either one or zero depending on whether $k$ divides $n$ or not. This
matrix is always triangular, since of course no divisor can exceed a given
number, and its determinant (for any dimension) is 1.%
\begin{equation}%
\begin{pmatrix}
\boldsymbol{n\diagdown k} & \mathbf{1} & \mathbf{2} & \mathbf{3} & \mathbf{4}
& \mathbf{5} & \mathbf{6} & \mathbf{7} & \mathbf{8} & \mathbf{9} & \mathbf{10}
& ...\\
\mathbf{1} & 1 & 0 & 0 & 0 & 0 & 0 & 0 & 0 & 0 & 0 & ...\\
\mathbf{2} & 1 & 1 & 0 & 0 & 0 & 0 & 0 & 0 & 0 & 0 & ...\\
\mathbf{3} & 1 & 0 & 1 & 0 & 0 & 0 & 0 & 0 & 0 & 0 & ...\\
\mathbf{4} & 1 & 1 & 0 & 1 & 0 & 0 & 0 & 0 & 0 & 0 & ...\\
\mathbf{5} & 1 & 0 & 0 & 0 & 1 & 0 & 0 & 0 & 0 & 0 & ...\\
\mathbf{6} & 1 & 1 & 1 & 0 & 0 & 1 & 0 & 0 & 0 & 0 & ...\\
\mathbf{7} & 1 & 0 & 0 & 0 & 0 & 0 & 1 & 0 & 0 & 0 & ...\\
\mathbf{8} & 1 & 1 & 0 & 1 & 0 & 0 & 0 & 1 & 0 & 0 & ...\\
\mathbf{9} & 1 & 0 & 1 & 0 & 0 & 0 & 0 & 0 & 1 & 0 & ...\\
\mathbf{10} & 1 & 1 & 0 & 0 & 1 & 0 & 0 & 0 & 0 & 1 & ...\\
... & ... & ... & ... & ... & ... & ... & ... & ... & ... & ... &
\end{pmatrix}
\label{Divisor matrix}%
\end{equation}
(Matrix (\ref{Divisor matrix}) is closely related to the Redheffer matrix, see
e.g. \cite{Redheffer}, \cite{Trott}.) Introducing%
\begin{equation}
\sigma_{0}(n):=%
{\displaystyle\sum\limits_{k=1}^{n}}
\alpha_{nk}\label{sigma}%
\end{equation}
we see that $\sigma_{0}(n)$ just counts the number of all divisors of a given
$n$ including both unity and $n$ itself.

It is known (see e.g. \cite{Inverse sequence}) that the inverse of matrix
(\ref{Divisor matrix}) is:%
\begin{equation}
\beta_{nk}=\left\{
\genfrac{}{}{0pt}{}{\mu\left(  \frac{n}{k}\right)  \quad\text{if}\quad k\mid
n}{0\quad\text{if}\quad k\nmid n}%
\right. \label{Beta}%
\end{equation}
where $\mu$ denotes the M\"{o}bius function:%
\begin{equation}
\mu(n)=\left\{
\begin{tabular}
[c]{l}%
$0$ if $n$ has squared prime factor\\
$+1$ if $n$ is a square-free positive integer with an even number of prime
factors\\
$-1$ if $n$ is a square-free positive integer with an odd number of prime
factors
\end{tabular}
\right. \label{Moebius mu}%
\end{equation}%
\begin{equation}%
\begin{pmatrix}
\boldsymbol{n\diagdown k} & \mathbf{1} & \mathbf{2} & \mathbf{3} & \mathbf{4}
& \mathbf{5} & \mathbf{6} & \mathbf{7} & \mathbf{8} & \mathbf{9} & \mathbf{10}
& ...\\
\mathbf{1} & 1 & 0 & 0 & 0 & 0 & 0 & 0 & 0 & 0 & 0 & ...\\
\mathbf{2} & -1 & 1 & 0 & 0 & 0 & 0 & 0 & 0 & 0 & 0 & ...\\
\mathbf{3} & -1 & 0 & 1 & 0 & 0 & 0 & 0 & 0 & 0 & 0 & ...\\
\mathbf{4} & 0 & -1 & 0 & 1 & 0 & 0 & 0 & 0 & 0 & 0 & ...\\
\mathbf{5} & -1 & 0 & 0 & 0 & 1 & 0 & 0 & 0 & 0 & 0 & ...\\
\mathbf{6} & 1 & -1 & -1 & 0 & 0 & 1 & 0 & 0 & 0 & 0 & ...\\
\mathbf{7} & -1 & 0 & 0 & 0 & 0 & 0 & 1 & 0 & 0 & 0 & ...\\
\mathbf{8} & 0 & 0 & 0 & -1 & 0 & 0 & 0 & 1 & 0 & 0 & ...\\
\mathbf{9} & 0 & 0 & -1 & 0 & 0 & 0 & 0 & 0 & 1 & 0 & ...\\
\mathbf{10} & 1 & -1 & 0 & 0 & -1 & 0 & 0 & 0 & 0 & 1 & ...\\
... & ... & ... & ... & ... & ... & ... & ... & ... & ... & ... &
\end{pmatrix}
\label{Inverse matrix}%
\end{equation}
Note that the numbers in (\ref{Inverse matrix}) when summed in rows give zero
except for the first row which stems from the following identity:%
\begin{equation}
\sum_{d|n}\mu(d)=\delta_{n,1}\label{Identity}%
\end{equation}
Matrices (\ref{Divisor matrix}) and (\ref{Inverse matrix}) are visualized in
Figure 1.

Somewhat similar but purely qualitative results have been published in
\cite{Cox}.%

\begin{center}
\includegraphics
{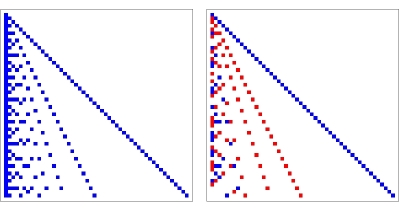}%
\\
Figure 1. Graphic distribution of divisors (\ref{Divisor matrix}) for
$n=1,2,...,50$ as a square matrix (left panel). Each blue square denotes $+1$.
In the inverse matrix (\ref{Inverse matrix}) blue square denotes $+1$ and red
square denotes $-1$ (right panel).
\end{center}
\bigskip

\section{General case}

The particular choice of the exponential function in (\ref{alpha3}) is not
crucial to our reasoning. Indeed, instead of this function we can take
\textit{any} regular function $f\left(  x\right)  $ provided that it has all
non-zero coefficients in its power series expansion%
\begin{equation}
f\left(  x\right)  =c_{0}+c_{1}x+c_{2}x^{2}+c_{3}x^{3}...\quad\quad c_{i}%
\neq0\text{ for }i=1,2,3,...\label{expansion}%
\end{equation}
Thus in general we have (up to appropriate normalizing factor)%
\begin{equation}
\alpha_{nk}=\left.  \frac{\mathrm{d}^{n}}{\mathrm{d}x^{n}}f\left(
x^{k}\right)  \right\vert _{x=0}\label{alpha8}%
\end{equation}

For example, taking%
\begin{equation}
f\left(  x\right)  =\frac{x}{1-x}=x+x^{2}+x^{3}+...\label{geometric}%
\end{equation}
we get:%
\begin{equation}
\alpha_{nk}=\left.
{\displaystyle\sum\limits_{j=1}^{\infty}}
\binom{jk}{n}x^{jk-n}\right\vert _{x=0}\label{alpha9}%
\end{equation}
or simply%
\begin{equation}
\alpha_{nk}=\left.
{\displaystyle\sum\limits_{j=1}^{\infty}}
\frac{x^{jk-n}}{\left(  jk-n\right)  !}\right\vert _{x=0}\label{alpha10}%
\end{equation}
Taking%
\begin{equation}
f\left(  x\right)  =\log\left(  1-x\right)  =-x-\frac{x^{2}}{2}-\frac{x^{3}%
}{3}-...\label{logarithm}%
\end{equation}
we get:%
\begin{equation}
\alpha_{nk}=\left(  -1\right)  ^{n/k}\frac{n}{k}\left.
{\displaystyle\sum\limits_{j=1}^{\infty}}
\frac{\left(  -1\right)  ^{j}}{j}\binom{jk}{n}x^{jk-n}\right\vert
_{x=0}\label{alpha11}%
\end{equation}

The general explicit formula for $\alpha_{nk}$ using arbitrary function $f$
satisfying (\ref{expansion}) is:%
\begin{equation}
\alpha_{nk}=\left(  \frac{n}{k}\right)  !\frac{1}{f^{(n/k)}\left(  0\right)
}\left.
{\displaystyle\sum\limits_{j=0}^{\infty}}
\frac{f^{(j)}\left(  0\right)  }{j!}\binom{jk}{n}x^{jk-n}\right\vert
_{x=0}\label{alpha12}%
\end{equation}
where $f^{(j)}\left(  0\right)  $ denotes the $j$-th derivative of $f$ with
respect to $x$ taken at $x=0$. (If $n/k$ in (\ref{alpha12}) is non-integer
then the value of fractional derivative $f^{(n/k)}\left(  0\right)  $\ is
unimportant since in this case the sum vanishes.)

The table below contains normalizing factors for $\alpha_{nk}$, for several
different choices of function $f(x)$, obtained using (\ref{alpha12}%
).\bigskip\bigskip%

\begin{tabular}
[t]{ll}%
$f\left(  x\right)  =e^{x}$ & $\alpha_{nk}=\frac{n}{k}!\frac{1}{n!}\left.
\frac{\mathrm{d}^{n}}{\mathrm{d}x^{n}}f\left(  x^{k}\right)  \right\vert
_{x=0}\bigskip$\\
$f\left(  x\right)  =\ln\left(  1-x\right)  $ & $\alpha_{nk}=-\frac{n}{k}%
\frac{1}{n!}\left.  \frac{\mathrm{d}^{n}}{\mathrm{d}x^{n}}f\left(
x^{k}\right)  \right\vert _{x=0}\bigskip$\\
$f\left(  x\right)  =\frac{x}{1-x}$ & $\alpha_{nk}=\frac{1}{n!}\left.
\frac{\mathrm{d}^{n}}{\mathrm{d}x^{n}}f\left(  x^{k}\right)  \right\vert
_{x=0}\bigskip$\\
$f\left(  x\right)  =\sqrt{1+x}$ & $\alpha_{nk}=-\frac{(-2)^{n/k}}{\left(
2\frac{n}{k}-3\right)  !!}\frac{n}{k}!\frac{1}{n!}\left.  \frac{\mathrm{d}%
^{n}}{\mathrm{d}x^{n}}f\left(  x^{k}\right)  \right\vert _{x=0}\bigskip$\\
$f\left(  x\right)  =\frac{1}{\sqrt{1+x}}$ & $\alpha_{nk}=\left(  -1\right)
^{n/k}\frac{\Gamma\left(  \frac{1}{2}\right)  }{\Gamma\left(  \frac{n}%
{k}+\frac{1}{2}\right)  }\frac{n}{k}!\frac{1}{n!}\left.  \frac{\mathrm{d}^{n}%
}{\mathrm{d}x^{n}}f\left(  x^{k}\right)  \right\vert _{x=0}\bigskip$\\
$f\left(  x\right)  =\left(  1+x\right)  ^{-3/2}$ & $\alpha_{nk}%
=\frac{(-2)^{n/k}}{\left(  2\frac{n}{k}+1\right)  !!}\frac{n}{k}!\frac{1}%
{n!}\left.  \frac{\mathrm{d}^{n}}{\mathrm{d}x^{n}}f\left(  x^{k}\right)
\right\vert _{x=0}\bigskip$\\
$f(x)=W(x)$ & $\alpha_{nk}=\left(  -1\right)  ^{n/k-1}\frac{\frac{n}{k}%
!}{\left(  \frac{n}{k}\right)  ^{\frac{n}{k}-1}}\frac{1}{n!}\left.
\frac{\mathrm{d}^{n}}{\mathrm{d}x^{n}}f\left(  x^{k}\right)  \right\vert
_{x=0}\bigskip$\\
$f\left(  x\right)  =\frac{1}{1-x-x^{2}}$ & $\alpha_{nk}=\frac{1}{F_{\frac
{n}{k}+1}}\frac{1}{n!}\left.  \frac{\mathrm{d}^{n}}{\mathrm{d}x^{n}}f\left(
x^{k}\right)  \right\vert _{x=0}\bigskip$%
\end{tabular}
\bigskip

($W(x)$ is the Lambert $W$-function and$\ F_{n}$ in the last row denotes the
$n$-th Fibonacci number.)

\section{Interpretation}

Let us now explain in more details how it all works. The thing is that all
formulas for $\alpha_{nk}$\ presented so far encode, at least to some extent,
the ancient algorithm known as the sieve of Eratosthenes.

Indeed, consider as $f(x)$\ the function $f(x)=x/(1-x)$ and let us temporarily
restrict ourselves to the linear case:$\ f(x)\approx x$. According to the
general formula (\ref{alpha8}) we have%
\begin{equation}
\alpha_{nk}=\frac{1}{n!}\left.  \frac{\mathrm{d}^{n}}{\mathrm{d}x^{n}}f\left(
x^{k}\right)  \right\vert _{x=0}=\frac{1}{n!}\left.  \frac{\mathrm{d}^{n}%
}{\mathrm{d}x^{n}}x^{k}\right\vert _{x=0}=\binom{k}{n}\left.  x^{k-n}%
\right\vert _{x=0}=\delta_{k,n}\label{alpha13}%
\end{equation}
and this produces a single line of ones on the diagonal $n=k$\ in the divisor
matrix (\ref{Divisor matrix}) -- cf. Figure 2 below. This is equivalent to the
trivial statement that all integers are divisible both by one and by
themselves. Let us further consider more precise approximation $f(x)\approx
x+x^{2}$. We get from (\ref{alpha8}) another sequence of ones on the line
$n=2k$. This is equivalent to selecting all even integers $n$ and adding to
the divisor matrix (\ref{Divisor matrix}) their divisors $n/2$. Taking into
account higher powers of $x$ we select all numbers $n$\ which are multiplies
of $3,4,5,...$ and this adds to the matrix further lines of divisors:
$n/3,n/4,n/5$, respectively.

Proceeding in the same way we finally arrive at the full expansion of $f(x)$:%
\begin{equation}
f(x)=\frac{x}{1-x}=%
{\displaystyle\sum\limits_{j=1}^{\infty}}
x^{j}\label{f}%
\end{equation}
which produces the entire sequence of lines $n=jk$ labelled by parameter
$j=1,2,3,...$. In this way we have selected and visualized \textit{all
divisors for all integers}. It is clear that there are certain well-defined
numbers $n$ (marked in bold in Figure 2) which have exactly two divisors:
unity and themselves, i.e. prime numbers: $2,3,5,7,11,13,...$ At the same time
we see the importance of condition $c_{i}\neq0$\ in (\ref{expansion}) since
even a single coefficient $c_{i}=0$ would cause a skipping of certain
divisors. In view of this the characteristic function for divisors may also be
written in a very natural form as a sum over Kronecker deltas:%
\begin{equation}
\alpha_{nk}=%
{\displaystyle\sum\limits_{j=1}^{n}}
\delta_{jk,n}\label{alpha14}%
\end{equation}

Note that combining (\ref{sigma}), (\ref{alpha8}) and (\ref{f}) gives:%
\begin{equation}
\sigma_{0}(n):=%
{\displaystyle\sum\limits_{k=1}^{n}}
\alpha_{nk}=\frac{1}{n!}\frac{\mathrm{d}^{n}}{\mathrm{d}x^{n}}\left.
{\displaystyle\sum\limits_{k=1}^{\infty}}
\frac{x^{k}}{1-x^{k}}\right\vert _{x=0}\label{Lambert1}%
\end{equation}
Hence%
\begin{equation}%
{\displaystyle\sum\limits_{k=1}^{\infty}}
\frac{x^{k}}{1-x^{k}}=%
{\displaystyle\sum\limits_{n=1}^{\infty}}
\sigma_{0}(n)x^{n}\label{Lambert2}%
\end{equation}
which is consistent with the theory of Lambert series (see e.g. \cite{Apostol}%
) which is the generating function for the sequence $\sigma_{0}(n)$ where
$\sigma_{0}(n)$\ is the total number of divisors for a given integer $n$.%

\begin{center}
\includegraphics
{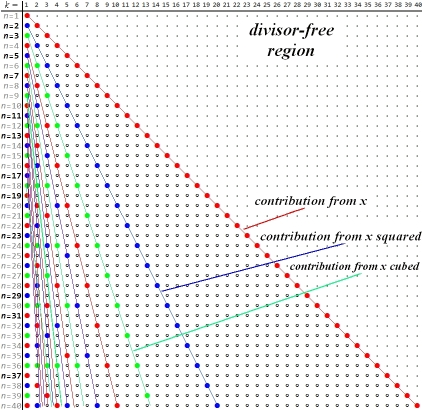}%
\\
Figure 2. Distribution of divisors of integers computed from $\alpha_{nk}$.
This figure illustrates how various terms in the sum (\ref{alpha14})
contribute to the whole pattern of divisors. Each term corresponds to a ray of
divisors.\ Rows are labelled by consecutive integers $n$ and columns are
labelled by potential divisors $k$. Each colored disc means that given $k$
actually divides $n$, otherwise there is small black circle. To better
visualize the whole pattern discs are in 3 different colors and lines
connecting them are drawn. Of course, above the diagonal ($k>n$) there can't
be any divisors.
\end{center}

\section{Concluding remarks}

A few elementary comments at the end of this note. As we have seen, all
divisors $k$ of integers $n$ lie on rays passing through the origin of the
coordinate system on the ($n,k$) plane and are labelled by an integer
parameter $j=1,2,3...$%
\begin{equation}
n=jk\label{condition}%
\end{equation}
We have also seen that this simple condition has a natural interpretation
since $j$ may be identified with the exponent in $x^{j}$ in the expansion
(\ref{expansion}).\ The key thing is that these rays must pass through certain
points of an integer lattice and only then a potential divisor\ can be an
actual divisor. For large $n$ these rays typically get closer and closer to
one another. Therefore we qualitatively see why it is so difficult to
factorize large integers.

Moreover, numerical experiments suggest that \textit{all} divisors lie on
countable families of parabolas passing through the origin (see Figures 3, 4
and 5 below). These parabolas are "quantized" in the sense that each family is
characterized by two discrete parameters $\mu=1,2,3...$ and $\nu=1,2,3...$ and
inside any family parabolas are labelled by another integer parameter $i$:%
\begin{equation}
g_{i}^{(\mu\nu)}(k)=-\frac{\mu}{\nu}k^{2}+\frac{i}{\nu}k\label{parabolas}%
\end{equation}
Careful simulations using Mathematica revealed that parameter $i$ assumes
equidistant values with integer constant step:%
\begin{equation}
\delta=\gcd(\mu,\nu)\label{step}%
\end{equation}
starting from $i=\mu+\nu$ where $\gcd$ denotes greatest common divisor, i.e.
$i=\mu+\nu$, $\mu+\nu+\delta$, $\mu+\nu+2\delta,...$%

\begin{center}
\includegraphics
{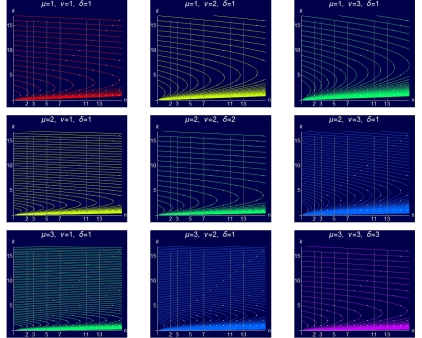}%
\\
Figure 3. Various families of parabolas (\ref{parabolas}) for $\mu=1,2,3$ and
$\nu=1,2,3$. Step $\delta$\ (\ref{step}) described in the main text is also
indicated.
\end{center}
%

\begin{center}
\includegraphics
{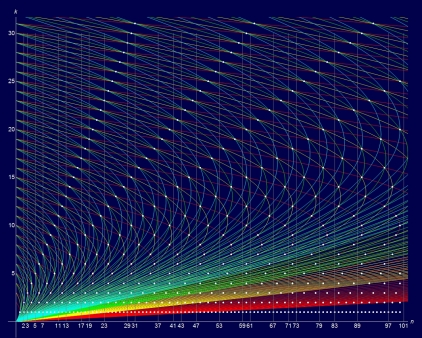}%
\\
Figure 4. Family of parabolas (\ref{parabolas}) for $\mu=1$ and $\nu=1$
(red)$,2$ (orange)$,3$ (green)$,$and $4$ (cyan) for $n<100$. For clarity of
the plot parameter $i$ assumes only $50$ consecutive\ values. Prime numbers
among $n$s are indicated by vertical lines.
\end{center}
%

\begin{center}
\includegraphics
{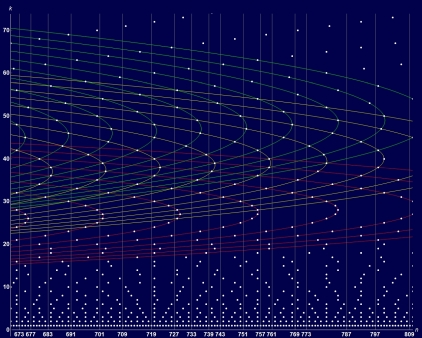}%
\\
Figure 5. Family of parabolas (\ref{parabolas}) for $\mu=1$ and $\nu=1$
(red)$,2$ (yellow) and $3$ (green) around $n=740$. For clarity of the plot
parameter $i$ assumes only $5$ consecutive\ values. Prime numbers among $n$s
are indicated by vertical lines.
\end{center}

As far as I am aware the unexpected parabolas in the distribution of divisors
have been independently noticed by Jeffrey Ventrella (see his popular book
\cite{Ventrella}, page 33) but with no quantitative considerations.

Finally, it should be stressed that, unfortunately, expressions presented in
this note do not tell us much about distribution of primes. They are even not
very suitable for numerical calculations for large $n$ therefore may be
treated merely as a curiosity. Nevertheless, we have shown some unexpected
relationship between number theory and calculus.\bigskip

\textbf{Acknowledgments.} \textit{The author would like to thank Prof. Jeffrey
Lagarias for his encouragement and several suggestions and to Prof. Andrzej
Schinzel for several remarks.}

\textit{The results presented in this paper were inspired by experimenting
with Wolfram Mathematica. Also all calculations were checked using this
powerful software.}

\end{document}